\documentclass[10pt,a4paper]{article}

\usepackage{amsmath, amsfonts, amsthm, amssymb, amscd}
\usepackage[mathcal]{euscript}
\usepackage{mathrsfs}

\newtheorem{theorem}{Theorem}[section]
\newtheorem{proposition}[theorem]{Proposition}

\newtheorem{corollary}[theorem]{Corollary}
\newtheorem{definition}[theorem]{Definition}
\newtheorem{remark}[theorem]{Remark}
\numberwithin{equation}{section}

\input{xypic.tex} \xyoption{all}
\usepackage{amsthm,amssymb,latexsym,amsmath}
\usepackage[all]{xypic}

\newcommand{\catA}[1]{{\mathfrak A}}
\newcommand{\catI}[1]{{\mathfrak I}}
\newcommand{\catS}[1]{{\mathfrak S}}

\newcommand{\im}{{\rm Im}~}

\DeclareMathOperator{\ho}{H}

\begin{document}

\title{
\centering  \large \bf
Monad constructions of omalous bundles
}

\author{  Abdelmoubine Amar Henni\footnote{Supported by FAPESP fellowship and grant 2009/18249-0}\and
 Marcos Jardim\footnote{Partially supported by the CNPq grant number 305464/2007-8 and FAPESP grant number 2005/04558-0.}
}

\date{}
\maketitle

\thispagestyle{empty}
\pagestyle{myheadings}
\markboth{A. A. Henni \quad M. Jardim}
  {Monad constructions of omalous bundles}
\setlength{\baselineskip}{6mm}

\begin{abstract}
We consider a particular class of holomorphic vector bundles relevant for supersymmetric string theory, called \emph{omalous}, over nonsingular projective varieties. We use monads to construct examples of such bundles over 3-fold hypersurfaces in $\mathbb{P}^{4}$, complete intersection Calabi-Yau manifolds in $\mathbb{P}^{k}$, blow-ups of $\mathbb{P}^{2}$ at $n$ distinct points, and products $\mathbb{P}^{m}\times\mathbb{P}^{n}$.
\end{abstract}

\section{Introduction}
\label{IN}

Let $X$ be a nonsingular projective variety, $TX$ be its tangent bundle and $\omega_{X}$ its canonical line bundle. This paper is dedicated to the study of the following class of holomorphic vector bundles.

\begin{definition}
A holomorphic vector bundle $\mathcal{E}\to X$ is called \emph{omalous} if it satisfies the following conditions:
\begin{itemize}
\item $c_{2}(\mathcal{E})=c_{2}(TX)$
\item $det(\mathcal{E}^{\ast})\simeq\omega_{X}$
\end{itemize}
\end{definition}

Recall also that a holomorphic vector bundle $\mathcal{E}\to X$ is \emph{slope (semi-)stable} with respect to a chosen polarization $\mathcal{O}_{X}(1)$ of $X,$ i.e., for every proper subsheaf $\mathcal{F}$ of $\mathcal{E}$ the inequality $\frac{deg(\mathcal{F})}{rk\mathcal{F}} (\leq)<\frac{deg(\mathcal{E})}{rk\mathcal{E}}$ is satisfied.

The nomenclature comes from the fact that the matching of the first and second Chern classes of $\mathcal{E}$ and $TX$ is the usual Green-Schwarz anomaly cancellation condition in heterotic string theory \cite{CHSW,GSW}; hence such bundles are \emph{not anomalous}, that is \emph{omalous} (Josh Guffin attributes this terminology to Ron Donagi, see the footnote in the first page of \cite{G}).

Such bundles have a long history in the string literature. They appeared in the attempt at compactifying superstring theory to a theory on a $M^{4}\times X,$ (where $X$ is complex compact Calabi-Yau 3-fold and $M^{4}$ is a flat Minkowsi space) with an unbroken $N=1$ supersymmetry in four dimensions. This was done by arguing pertubatively in terms of the string coupling constant as in \cite{CHSW, Wittens} and \cite{GSW}.

In the late nineties, arguments about compactifications using non-perturbative vacua in heterotic string theory were given by introducing five-branes as done by Donagi et. al. in \cite{Donagi}. In this context, the general formula for anomaly cancelation is given by
$$ c_{2}(\mathcal{E})-c_{2}(TX)=[W], $$
where $[W]$ is the cohomology class of a four-form on the five-brane, and by Poincar\'e duality it corresponds to the class of an effective curve in the five-brane \cite[Section 2]{Donagi}.

More recently, the motivation to consider such bundles comes from two sources. First, the omality conditions are necessary to the construction of a \emph{quantum sheaf cohomology} for the bundle $\mathcal{E}$, c.f. \cite[Section 2]{G} and \cite{GK}. The quantum sheaf cohomology of a bundle $\mathcal{E}\to X$ is a generalization of the quantum cohomology of $X$, and consists of the structure of a Frobenius algebra on
$$ Q\ho^\bullet(\mathcal{E}):= \oplus_{p,k} \ho^p(X,\wedge^k\mathcal{E}^\ast) \otimes \mathbb{C}[[\underline{q}]] $$
with product and bilinear pairing induced by the three-point correlation functions in a $(0,2)$ supersymmetric nonlinear sigma model; see \cite{G,GK} and the references therein for further details.

Second, it was shown by Andreas and Garcia-Fernandez in \cite{AGF} that stable omalous bundles over Calabi-Yau 3-folds admit solutions of the \emph{Strominger system}, which is a system of coupled partial differential equations defined over a compact complex manifold relevant in heterotic string theory, c.f. \cite{AGF} and the references therein.

The simplest example of omalous bundles are $TX\oplus\mathcal{O}_{X}^{\oplus k}$ and its deformations; a few other examples were considered in \cite{AGF,B,DZ,GK}. Our goal is to construct more examples of (stable) omalous bundles over various choices for $X$ using \emph{monads}.

We remark that in the attempt of giving phenomenological models from heterotic compactifications, particular \emph{monads} have been used in the physics literature such as in \cite{Anderson1, Anderson2, Anderson3, Blumen, Kachru}; we emphasize however that what is called a monad in \cite{Anderson1, Anderson2, Anderson3} does not coincide with the usual definition in the mathematical literature.

Recall that a \emph{monad} on $X$ is a complex of locally free sheaves
$$ {\rm M}_\bullet ~~:~~ M_{0} \stackrel{\alpha}{\longrightarrow} M_1 \stackrel{\beta}{\longrightarrow} M_2 $$
such that $\beta$ is locally right-invertible, $\alpha$ is locally left-invertible. The (locally free) sheaves $K:=ker\beta,$ $Q:=coker\alpha$ and $E:=ker\beta/\im\alpha$ are called, respectively, the \emph{kernel}, \emph{cokernel} and \emph{cohomology} of ${\rm M}_\bullet$.

In what follows, we provide examples of omalous bundles over 3-fold hypersurfaces in $\mathbb{P}^{4}$, complete intersection Calabi-Yau manifolds in $\mathbb{P}^{k}$ ($k=4,5,6,7$), blow-ups of $\mathbb{P}^{2}$ at $n$ distinct points, and products $\mathbb{P}^{m}\times\mathbb{P}^{n}$. All of these examples arise as cohomology, kernel, or cokernel of particular monads over these manifolds. We hope that such examples will be relevant for a deeper understanding of both quantum sheaf cohomology, the Strominger system and supersymmetric string theory.

\paragraph{Acknowledments}
This paper would not exist if we had not heard Josh Guffin's and Mario Garcia-Fernandez's excelent talks during the Second Latin Congress on Symmetries in Geometry and Physics; we thank them for useful conversations during the conference. We also thank Rosa Maria Mir\'o-Roig for her comments on the first version of this paper.

\section{Stable omalous bundles over 3-fold hypersurfaces in $\mathbb{P}^{4}$}

Let $X$ be the non-singular quintic $3-$fold in $\mathbb{P}^{4},$ and consider the following monad:
\begin{equation}\label{m-quintic}
\xymatrix{0\ar[r]&\mathcal{O}_{X}(-1)^{\oplus10}\ar[r]^{\alpha}&\mathcal{O}_{X}^{\oplus22}\ar[r]^{\beta}&\mathcal{O}_{X}(1)^{\oplus10}\ar[r]&0}.
\end{equation}
The existence of the monads (\ref{m-quintic}) is explicitly guaranteed by the construction in \cite[Section 3]{Jardim}; its cohomology bundle $\mathcal{E}$ is a rank $2$ bundle with Chern classes $c_{1}(\mathcal{E})=0$ and $c_{2}(\mathcal{E})=10\cdot H^{2},$ where $H=c_{1}(\mathcal{O}_{X}(1))$ is the ample generator of the Picard group of $X.$ This bundle is stable by \cite[Main Theorem]{Jardim}. Moreover $c_{2}(\mathcal{E})=c_{2}(TX)$ and $det(\mathcal{E}^{\ast})\cong\mathcal{O}_{X}=\omega_{X}$ since the quintic $3-$fold in $\mathbb{P}^{4}$ is Calaby-Yau. Hence $\mathcal{E}$ is a stable omalous bundle over $X.$

\bigskip

Let us consider now a non-singular $3-$fold $X_{d}$ of degree $d$ in $\mathbb{P}^{4}.$ One can show that $$c_{1}(TX_{d})=(5-d)\cdot H\textnormal{ and}$$ $$c_{2}(TX_{d})=(d^{2}-5d+10)\cdot H^{2}.$$
Let $\mathcal{E}$ be a rank $3$ linear bundle, that is, the cohomology of a linear monad of the form
\begin{equation}\label{linear}
\xymatrix{0\ar[r]&\mathcal{O}^{\oplus (c+l)}_{X_{d}}(-1)\ar[r]^{\alpha}&\mathcal{O}^{\oplus(3+2c+l)}_{X_{d}}\ar[r]^{\beta}&\mathcal{O}^{\oplus c}_{X_{d}}(1)\ar[r]&0}.
\end{equation}
Then $c_{1}(\mathcal{E})=l\cdot H$ and $c_{2}(\mathcal{E})=[\frac{l}{2}(l+1)+c]\cdot H^{2}.$
\begin{proposition}\label{Omalous solution}
\begin{itemize}
\item[(i)] The cohomology bundle $\mathcal{E}$ of the linear monad \eqref{linear} is omalous for every odd integer $k\geq7,$ such that $(d,l,c)$ are given by
$$d(k)=\frac{1}{2}(k-1),\quad\quad l(k)=\frac{1}{2}(11-k),\quad\quad c(k)=\frac{1}{8}(k^{2}-41).$$
\item[(ii)] Furthermore the omalous bundle $\mathcal{E}$ is stable for $(d,l,c)=(3,2,1),(4,1,5)$ and semi-stable for $(d,l,c)=(5,0,10).$
\end{itemize}
\end{proposition}

\vspace{.2in}

\noindent{\bf Proof.}
$(i)$ The conditions for which $\mathcal{E}$ is omalous are given by $c_{1}(\mathcal{E})=c_{1}(TX)$ and $c_{2}(\mathcal{E})=c_{2}(TX).$ In this case, one must have
$$
\left\{\begin{array}{l}d^{2}-5d+10=c+\frac{1}{2}(5-d)(6-d) \\ l=5-d \end{array}\right.
$$
Thus one must look for positive integer solutions $d(c)$ of the quadratic equation $d^{2}+d-(10+2c)=0,$ i.e., $d=\frac{1}{2}(-1+\sqrt{41+8c}).$ For every odd integer $k\geq7,$ it is easy to verify that $d(k)=\frac{1}{2}(k-1)$ are the desired roots, and the corresponding value for $c$ is $c(k)=\frac{1}{8}(k^{2}-41)$.

\bigskip

$(ii)$ The stability part follows from \cite[Theorem 7]{Jardim2}. If $d=5$ then $l=0$ and $c=10,$ Then $\mathcal{E}$ is an instanton bundle, and by \cite[Theorem 3]{Jardim2}, it is semi-stable.

\hfill $\Box$

\begin{remark}\rm
The existence of the monads (\ref{linear}) is a consequence of Flyostad's Theorem for monads on $\mathbb{P}^{4}$. More precisely, the Main Theorem of \cite{F} implies that the degeneration locus of a generic monad of the form
$$ \xymatrix{0\ar[r]&\mathcal{O}^{\oplus (c+l)}_{\mathbb{P}^{4}}(-1)\ar[r]^{\alpha}&\mathcal{O}^{\oplus(3+2c+l)}_{\mathbb{P}^{4}}\ar[r]^{\beta}&\mathcal{O}^{\oplus c}_{\mathbb{P}^{4}}(1)\ar[r]&0}. $$
is zero dimensional. Its restriction to a generic hypersurface $X_d$ is precisely (\ref{linear}), hence its cohomology yields a vector bundle over it.
\end{remark}

\section{Omalous bundles on complete intersection \\ Calabi-Yau $3-$folds}
Let $X$ be a complete intersection Calabi-Yau $3-$fold in $\mathbb{P}^{n}.$ There are only five such cases, namely:
\begin{itemize}
\item A quintic in $\mathbb{P}^{4}.$
\item In $\mathbb{P}^{5},$ either the intersection of two cubics or the intersection of a quadric and a quartic.
\item In $\mathbb{P}^{6}$ the intersection of two quadrics with a cubic.
\item In $\mathbb{P}^{7}$ the intersection of four quadrics.
\end{itemize}

One can write $X=\cap_{i}^{l}X_{i}$ where the $X_{i}$'s are given as in the list above and $l=codim_{\mathbb{P}^{n}}(X)$. Moreover one has the following short exact sequences $$0\longrightarrow TX_{i}\longrightarrow T\mathbb{P}^{n}|_{X_{i}}\longrightarrow\mathcal{N}_{i}\longrightarrow0$$ $$0\longrightarrow\mathcal{O}_{\mathbb{P}^{n}}(-d_{i})\longrightarrow\mathcal{O}_{\mathbb{P}^{n}}\longrightarrow\mathcal{O}|_{X_{i}}\longrightarrow0$$ where the normal bundle $\mathcal{N}_{i}$ to $X_{i}$ is simply the invertible sheaf $\mathcal{O}_{\mathbb{P}^{n}}(d_{i})$ since each of the $X_{i}$ is a hypersurface of degree $d_{i}$ in $\mathbb{P}^{n}$. Using these sequences one can easily prove that the Chern Class of the tangent bundle $TX,$ to $X,$ is given by the formula $$C(TX)=\frac{(1+h)^{n+1}}{\Pi_{i=1}^{l}(1+d_{i}h)}$$ where $h=c_{1}(\mathcal{O}_{\mathbb{P}^{n}}(1)).$ From the condition $c_{1}(TX)=0,$ it follows that $$c_{2}(TX)=\frac{1}{2}[(\Sigma_{i=1}^{l}d_{i}^{2})-(n+1)]h^{2}.$$

\begin{proposition}
Let $\mathcal{E}$ be a rank $2$ bundle on a complete intersection Calabi-Yau $3-$fold $X$ given by the cohomology of the following monad $$\xymatrix{M:&0\ar[r]&\mathcal{O}_{X}(-1)^{\oplus c}\ar[r]^{\alpha}&\mathcal{O}_{X}^{\oplus2+2c}\ar[r]^{\beta}&\mathcal{O}_{X}(1)^{\oplus c}\ar[r]&0}.$$ with $c$ given according to the following table:
\begin{center}
\begin{tabular}{|c|c|c|c|} \hline
$X$&l&$(d_{1},\cdots,d_{l})$& $c=\frac{1}{2}[(\Sigma_{i=1}^{l}d_{i}^{2})-(n+1)]$ \\ \hline
$X\subset\mathbb{P}^{4}$&$1$&$5$&$10$ \\ \hline
$X\subset\mathbb{P}^{5}$&$2$&$ (3,3) $& $6 $ \\ \hline
$X\subset\mathbb{P}^{5}$&$2$&$ (4,2)$& $ 7$ \\ \hline
$X\subset\mathbb{P}^{6}$&$3$&$(2,2,3)$&$5$ \\ \hline
$X\subset\mathbb{P}^{7}$&$4$&$(2,2,2,2)$&$4$ \\ \hline
\end{tabular}
\end{center}
Then $\mathcal{E}$ is a stable and omalous.
\end{proposition}
\vspace{.2in}

\noindent{\bf Proof.}
Follows from the Main Theorem in \cite{Jardim} and the calculations above.
\hfill $\Box$

\section{Omalous bundles on multi-blow-ups of the projective plane}

Let $\pi:\tilde{\mathbb{P}}(n)\longrightarrow\mathbb{P}^{2}$ be the blow-up of the projective plane at $n$ distinct points. Its Picard group is generated by $n+1$ elements, namely: $Pic(\tilde{\mathbb{P}}(n))=\oplus_{i=1}^{n}E_{i}\mathbb{Z}\oplus H\mathbb{Z},$ where every $E_{i}$ is an exceptional divisor and $H$ is the divisor given by the pull-back of the generic line in $\mathbb{P}^{2}$. The intersection form is given by: $E_{i}^{2}=-1$, $E_{i}\cdot E_{j}=0$ for $i\neq j$, $E_{i}\cdot H=0$ and $H^{2}=1.$ The canonical divisor of the surface $\tilde{\mathbb{P}}(n)$ is given by $K_{\tilde{\mathbb{P}}(n)}=-3H+\Sigma_{i=1}^{n}E_{i}.$
In terms of line bundles, a divisor of the form $D=pH+\Sigma_{i=1}^{n}q_{i}E_{i}$ has the associated line bundle $\mathcal{O}(D)=\mathcal{O}(p,\overrightarrow{q})=\mathcal{O}(pH)\otimes\mathcal{O}(q_{1}E_{1})\otimes\cdots\otimes\mathcal{O}(q_{n}E_{n})$ where $\overrightarrow{q}=(q_{1},\cdots , q_{n}).$

Let $H^{2}\in \ho^{4}(\tilde{\mathbb{P}}(n),\mathbb{Z})$ be the fundamental class of $\tilde{\mathbb{P}}(n)$. For a torsion-free sheaf $\mathcal{E}$, of Chern character $ch(\mathcal{E})=r+(aH+\Sigma_{i=1}^{n}a_{i}E_{i})-(k-\frac{a^{2}-|\overrightarrow{a}|^{2}}{2})H^{2}$, twisted by a line bundle $\mathcal{O}(p,\overrightarrow{q})$ the Riemann-Roch formula is given by:
$$\chi(\mathcal{E}(p,\overrightarrow{q}))=-[k-\frac{a}{2}(a+3)+\frac{1}{2}\Sigma_{i=1}^{n}a_{i}(a_{i}-1)]+
\frac{r}{2}[(p+1)(p+2)-\Sigma_{i=1}^{n}q_{i}(q_{i}-1)]$$ $$+[ap-\Sigma_{i=1}^{n}a_{i}q_{i}].$$ Note that the notations used through this section are the ones given in \cite{henni}.
Omalous bundles $\mathcal{E}$ on $\tilde{\mathbb{P}}(n)$ are given in this case by the conditions: $$det(\mathcal{E}^{\ast})=\omega_{\tilde{\mathbb{P}}(n)}=\mathcal{O}_{\tilde{\mathbb{P}}(n)}(-3,\vec{1}),\quad\quad c_{2}(\mathcal{E})=c_{2}(T\tilde{\mathbb{P}}(n))=(3+n)\cdot H^{2}$$ Moreover, suppose that the direct image $\pi_{\ast}(\mathcal{E}),$ of $\mathcal{E},$ is a normalized and semi-stable torsion free sheaf (in our case, normalized means that $3<r$). Then we have the following:
\begin{proposition}
On a multi-blow-up $\tilde{\mathbb{P}}(n)$ of the projective plane with $n\geq3,$ let $\mathcal{E}$ be an omalous bundle of rank $r>3$ with semi-stable direct image $\pi_{\ast}(\mathcal{E}).$ Then $\mathcal{E}$ is the cohomology of the following monad
$$\xymatrix@C-0.5pc{M:& 0\ar[r]& \oplus_{i=0}^{n}K_{i}(-1,E_{i})\ar[r]^{\quad\alpha}& W\otimes\mathcal{O}_{\tilde{\mathbb{P}}(n)} \ar[r]^{\beta\quad}&\oplus_{i=0}^{n}L_{i}(1,-E_{i})\ar[r]&0}$$
where we put $E_{0}:=0$ and $K_{i}, L_{i}$ and $W$ are vector space of dimensions
$$\dim K_{i}=\left\{\begin{array}{ll}n & i=0\\ 2n-3 & \textnormal{otherwise}
\end{array}\right.\qquad\dim L_{i}=\left\{\begin{array}{ll}2n-3 & i=0\\ 2n-4 & \textnormal{otherwise}
\end{array}\right.$$ and $\dim W=4n(n-1)-3+r.$
\end{proposition}

\vspace{.2in}

\noindent{\bf Proof.}
The existence of the monad is guaranteed by \cite[Proposition 1.10]{Buch} since the direct image $\pi_{\ast}(\mathcal{E})$ is semi-stable and normalized ($r>3$). The omality condition implies that the bundle $\mathcal{E}$ has the following Chern character $ch(\mathcal{E})=r+(3H-\Sigma_{i=1}^{n}E_{i})-\frac{3}{2}(n-1)H^{2}.$ The vector spaces in the monad are explicitly given by \cite[Proposition1.10]{Buch}:
$K_{0}=\ho^{1}(\tilde{\mathbb{P}}(n),\mathcal{E}^{\ast}(-1,0)),$ $K_{i}=\ho^{1}(\tilde{\mathbb{P}}(n),\mathcal{E}(-1,0))$ for $i\neq0,$
and $L_{0}=\ho^{1}(\tilde{\mathbb{P}}(n),\mathcal{E}(-1,0)),$ $L_{i}=\ho^{1}(\tilde{\mathbb{P}}(n),\mathcal{E}(-1,E_{i}))$ for $i\neq0.$
Their dimensions follow by the Riemann-Roch Formula.

\hfill $\Box$

\section{Omalous bundles on $\mathbb{P}^{n}\times\mathbb{P}^{m}$}
Let $X=\mathbb{P}^{n}\times\mathbb{P}^{m}, $ with the natural projections $$\xymatrix@C-1pc{&\mathbb{P}^{n}\times\mathbb{P}^{m}\ar[ld]_{\pi_{1}}\ar[rd]^{\pi_{2}}& \\ \mathbb{P}^{n}& &\mathbb{P}^{m}}$$ Its Picard group is generated by $h_{1}=\pi_{2}^{\ast}c_{1}(\mathcal{O}_{\mathbb{P}^{n}}(1))$ and $h_{2}=\pi_{2}^{\ast}c_{1}(\mathcal{O}_{\mathbb{P}^{m}}(1)),$ then $Pic(X)=h_{1}\mathbb{Z}\oplus h_{2}\mathbb{Z}.$ The Chow ring of $X$ is given by $$A(X)=\mathbb{Z}[h_{1},h_{2}]\slash(h_{1}^{n+1}, h_{2}^{m+1}).$$ Let $L:=\mathcal{O}_{X}(1,1)$ be the ample line bundle associated to the ample divisor $h_{1}+h_{2}.$ For any sheaf $\mathcal{F}$ on $X,$ we define the degree of $\mathcal{F}$ with respect to $L$ as $deg_{L}(\mathcal{F}):=c_{1}(\mathcal{F})\cdot c_{1}(L)^{n+m-1}.$ Using the binomial formula one obtains $c_{1}(L)^{n+m-1}=l(n,m)[h_{1}^{n-1}\cdot h_{2}^{m}+\frac{m}{n}h_{1}^{n}\cdot h_{2}^{m-1}]$ where $l(n,m)=\frac{n(n+1)\cdots(n+m+1)}{m!}.$ It follows that if $c_{1}(\mathcal{F})=p\cdot h_{1}+q\cdot h_{2},$ then $deg_{L}(\mathcal{F})=l(n,m)[p+\frac{m}{n}q]\cdot h_{1}^{n}\cdot h_{2}^{m}.$ We define the $L-$\emph{slope} $\mu_{L}(\mathcal{F})$ of the sheaf $\mathcal{F}$ by $\mu_{L}(\mathcal{F}):=\frac{deg_{L}(\mathcal{F})}{rk(\mathcal{F})},$ and will say that $\mathcal{F}$ is $L-$\emph{(semi-)stable} if for every subsheaf $\mathcal{G}\subset\mathcal{F}$ the inequality $\mu_{L}(\mathcal{G})(\leq)<\mu_{L}(\mathcal{F})$ is satisfied.

%


The tangent bundle of $X$ is given by the following Euler sequence:
$$0\longrightarrow\mathcal{O}_{X}^{\oplus 2}\longrightarrow\mathcal{O}_{X}^{\oplus n+1}(1,0)\oplus\mathcal{O}_{X}^{\oplus m+1}(0,1)\longrightarrow TX\longrightarrow0,$$
from which one can easily compute the canonical bundle $\omega_{X}=\mathcal{O}_{X}(-n-1,-m-1),$ the first Chern class $c_{1}(TX)=c_{1}(X)=(n+1)\cdot h_{1}+(m+1)\cdot h_{2}$ and the second Chern class $c_{2}(TX)=\frac{1}{2}(n+1)h^{2}_{1}+\frac{1}{2}m(m+1)h^{2}_{2}+(n+1)(m+1)\cdot h_{1}\cdot h_{2}.$

Now let us consider a rank $(b+c-a)-$bundle $Q$ fitting in the following short exact sequence:
\begin{equation}\label{coker}
0\longrightarrow\mathcal{O}_{X}^{\oplus a}\longrightarrow\mathcal{O}_{X}(1,0)^{\oplus b}\oplus\mathcal{O}_{X}(0,1)^{\oplus c}\longrightarrow Q\longrightarrow0.
\end{equation}

\begin{proposition}\label{L-stable}
\begin{itemize}
\item[(i)] $Q$ is omalous for the values $(b,c)=(n+1,m+1).$
\item[(ii)] The bundle $Q$ is $L$-stable.
\end{itemize}
\end{proposition}

\vspace{.2in}

\noindent{\bf Proof.}
\begin{itemize}
\item[(i)] From the exact sequence defining the bundle $Q$ one can easily compute the Chern classes $c_{1}(Q)=b\cdot h_{1}+c\cdot h_{2}$ and $c_{2}(Q)=\frac{b(b-1)}{2}\cdot h_{1}^{2}+\frac{c(c-1)}{2}\cdot h_{1}+bc\cdot h_{1}\cdot h_{1}.$ The result follows by imposing the conditions $c_{1}(Q)=c_{1}(X)$ and $c_{2}(Q)=c_{2}(TX).$

\item[(ii)] Note that the twisted dual bundle $Q^{\ast}(0,1)$ fits in the following exact sequence
\begin{equation}\label{m-pnpm}
0\longrightarrow Q^{\ast}(0,1)\longrightarrow\mathcal{O}_{X}^{\oplus b}(-1,1)\oplus\mathcal{O}_{X}^{\oplus c}\longrightarrow\mathcal{O}_{X}(0,1)^{\oplus a}\longrightarrow0.
\end{equation}
By using the first statement in \cite[Theorem 8]{JSE}, it follows that $Q^{\ast}(0,1)$ is $L$-stable, thus $Q$ is also $L$-stable.
\end{itemize}
\hfill $\Box$

In particular, when $a=2$ and the omalous conditions are satisfied, then $Q$ has exactly rank $n+m,$ thus it is a deformation of the tangent bundle $TX.$ Moreover it is easy to see that any deformation of $TX$ is given as the last term bundle in the sequence \eqref{coker}, with $a=2,$ $b=n+1$ and $c=m+1.$ Hence we have the following:
\begin{corollary}
Any deformation of the tangent bundle $TX$ of $X=\mathbb{P}^{n}\times\mathbb{P}^{m}$ is $L$-stable.
\end{corollary}

\vspace{.2in}

\noindent{\bf Proof.}
Follows from $(i)$ in \ref{L-stable}.

\hfill $\Box$


\vspace{.4in}


\vspace{.4in}

\begin{tabbing}
---------------------------\=-----------------------\=------------------\kill
Abdelmoubine Amar Henni     \>      \> Marcos Jardim \\
IMECC - UNICAMP             \>      \> IMECC - UNICAMP  \\
Departamento de Matematica  \>      \> Departamento de Matematica \\
Rua Sergio Buarque de Holanda, 651\>       \> Rua Sergio Buarque de Holanda, 651 \\
13083-859, Campinas-SP,  \>    \> 13083-859, Campinas-SP,   \\
Brazil \>          \>  Brazil  \\
e-mail: henni@ime.unicamp.br  \>       \> e-mail: jardim@ime.unicamp.br
\end{tabbing}


\vspace{.3in}

\begin{thebibliography}{99}

\bibitem{AGF}
Andreas, B., Garcia-Fernandez, M.,
{\em Solutions of the Strominger system via stable bundles on Calabi-Yau
threefolds},
Preprint arXiv: 1008:1018.

\bibitem{Anderson1}
Anderson, L. B., He, Y. H., Lukas, A.,
{\em Heterotic compactification, an algorithm approach},
JHEP 0707:049
J. High Energy Phys. 2007, no. 7, 049.

\bibitem{Anderson2}
Anderson, L. B., He, Y. H., Lukas, A.,
{\em Monad bundles in heterotic string compactification},
JHEP 0807:104
J. High Energy Phys. 2008, no. 7, 104.

\bibitem{Anderson3}
Anderson, L. B., Gray, J., He, Y. H., Lukas, A.,
{\em Exploring positive monad bundles and a new heterotic standard model},
JHEP 1002:054
J. High Energy Phys. 2010, no. 2, 054.

\bibitem{Blumen}
Blumenhagen, R.,
{\em Target space duality for $(0,2)$ compactifications},
Nucl. Phys. B 513, 573 (1998).

\bibitem{B}
Brambilla, M. C.,
{\em Semistability of certain bundles on a quintic Calabi-Yau threefold},
Rev. Mat. Complut. {\bf 22} (2009), 53--61.

\bibitem{Buch}
Buchdahl, N. P.,
{\em Monads and bundles on rational surfaces},
Rocky Mountain J. Math. {\bf 34} (2004), 513--540.

\bibitem{CHSW}
Candelas, P., Horowitz, G. T., Strominger, A., Witten, E.,
{\em Vacuum configurations for superstrings},
Nucl. Phys. B 258, 46 (1985).

\bibitem{Donagi}
Donagi, R., Lukas, A., Ovrut, B. A., Waldram, D.,
{\em Holomorphic vector bundles and non-perturbative vacua in M-theory},
JHEP 9906:034,, 1999.
J. High Energy Phys. 1999, no. 6, 034.

\bibitem{DZ}
Douglas, M. R., Zhou, C.-G.,
{\em Chirality Change in String Theory},
JHEP 06:014, 2004.
J. High Energy Phys. 2004, no. 6, 014.

\bibitem{F}
Floystad, G.,
Monads on projective spaces.
Comm. Algebra {\bf 28} (2000), 5503-5516.

\bibitem{GSW}
Green, M. B., Schwarz, J. H., Witten, E.,
{\em Superstring theory Vol.2: Loop amplitudes, anomalies and phenomenology}
Cambridge Monographs on Mathematical Physics. Cambridge University press, (1987)

\bibitem{G}
Guffin, J.,
{\em Quantum Sheaf Cohomology, a pr\'ecis},
Preprint arXiv:1101.1305.

\bibitem{GK}
Guffin, J., Katz, S.,
{\em Deformed quantum cohomology and (0,2) mirror symmetry},
J. High Energy Phys. 2010, no. 8, 109.

\bibitem{henni}
Henni, A. A.,
{\em Monads for torsion-free sheaves on multi-blow-ups of the projective plane},
Preprint arXiv:0903.3190.

\bibitem{Jardim}
Jardim, M.,
{\em Stable bundles on 3-fold hypersurfaces},
Bull. Braz. Math. Soc,  {\bf 38} (2007), 649--659.

\bibitem{Jardim2}
Jardim, M.,   Mir\'o-Roig, R. M.,
{\em On the stability of instanton sheaves over certain projective varieties},
Comm. Algebra {\bf 36} (2008), 288--298.

\bibitem{JSE}
Jardim, M.,  S\'a Earp, H. N.,
{\em Monad constructions of asymptotically stable bundles},
To appear.

\bibitem{Kachru}
Kachru, S.,
{\em Some three generation $(0,2)$ Calabi-Yau models},
Phys. Lett. B 349, 76 (1995).

\bibitem{Wittens}
Witten, E., Witten, L.,
{\em Large radius expension of superstring compactification},
Nucl. Phys. B 281, 109 (1987).










\end{thebibliography}
\end{document}